\documentclass[11pt,reqno]{amsproc} 
\usepackage{amsmath}
\usepackage{amssymb} 
\usepackage[all]{xy} 
\usepackage{psfrag}
\usepackage{graphicx}

\swapnumbers 
\theoremstyle{plain}
\newtheorem{theorem}{Theorem}[section]
\newtheorem{prop}[theorem]{Proposition}
\newtheorem{lemma}[theorem]{Lemma}
\newtheorem{cor}[theorem]{Corollary}

\theoremstyle{definition}
\newtheorem*{ack}{Acknowledgements}
\newtheorem{defn}[theorem]{Definition}

\newtheorem{example}[theorem]{Example}
\newtheorem{remark}[theorem]{Remark}
\newtheorem{notation}[theorem]{Notation} 

\newtheoremstyle{citing}
  {}
  {}
  {\itshape}
  {\parindent}
  {\scshape}
  {.}
  { }
  {\thmnote{#3}}

\theoremstyle{citing}

\newcommand{\norm}[2]{\left\| {#1}\right\| _{#2}}
\newcommand{\abs}[1]{\left |{#1}\right |}
\newcommand{\cc}{\subset\!\!\!\subset}
\renewcommand{\l}{{\mathcal L}} 
\newcommand{\ls}[1]{{\mathcal L}_{{#1}}} 
\newcommand{\set}[1]{\left\{{#1}\right \}}

\newcommand{\all}[2]{\left \{ {#1}\,|\,{#2} \right \}}

\newcommand{\I}{\mathcal I}
\newcommand{\C}{\mathbb C}

\newcommand{\N}{\mathbb N} 

\newcommand{\D}{\mathcal D}

\newcommand{\un}{\underline{n}}

\newcommand{\bt}{t}

\title[On the Ruelle eigenvalue sequence]{On
the Ruelle eigenvalue sequence}

\author{Oscar F.~Bandtlow and Oliver Jenkinson}

\address{Oscar Bandtlow; 
School of Mathematical Sciences, Queen
  Mary, University of London, Mile End Road, London, E1 4NS, UK.
\newline
{\tt ob@maths.qmul.ac.uk} \newline {\tt
    www.maths.qmul.ac.uk/$\sim$ob}}  

\address{Oliver Jenkinson; 
School of Mathematical Sciences, Queen
  Mary, University of London, Mile End Road, London, E1 4NS, UK.
  \newline {\tt omj@maths.qmul.ac.uk} \newline {\tt
    www.maths.qmul.ac.uk/$\sim$omj}}

\begin{document}

\date{\today}

\begin{abstract}
For certain real analytic data, we show that the
eigenvalue sequence 
of the associated transfer operator $\l$ is 
insensitive to the holomorphic function space on
which $\l$ acts.
Explicit bounds on this eigenvalue
sequence are established.
\end{abstract}

\maketitle

\section{Introduction}
\label{introsection}

For compact $X\subset\C^d$, and appropriate
real analytic $T_i:X\to X$ and $w_i:X\to\C$,
Ruelle \cite{ruelleinventiones}
considered the action of the transfer operator 
$\l f:= \sum_{i} w_i\cdot f\circ T_i$ on $U(D)$,
where $D\subset\C^d$ is a domain on which all the 
$T_i$'s and $w_i$'s are holomorphic, and $U(D)$  consists
of those holomorphic functions on $D$
which extend continuously to the closure of $D$.
Ruelle proved that $\l:U(D)\to U(D)$ 
is nuclear, hence in particular compact,
and that its
eigenvalue sequence
$\{\lambda_n(\l)\}_{n=1}^\infty$, henceforth referred to as
the \emph{Ruelle eigenvalue sequence}, is given by
the reciprocals
of the zeros of a dynamical determinant $\Delta$
(see (\ref{ruelledet}) for the definition).

In view of its various interpretations and applications
(e.g.~correlation decay rates \cite{baladibook,cpr},
Fourier resonances \cite{ruelleresonances2}, Laplacians
for hyperbolic surfaces 
\cite{pollicottadvmath, pollicottrocha}, Feigenbaum
period-doubling \cite{aac,ccr,jms,pollicott}),
it is desirable to establish explicit bounds on the
Ruelle eigenvalue sequence.
In the case where 
$D$ may be chosen as a ball,
and the $T_i$ all
map $D$ within the concentric ball
whose radius is $r<1$ times that of $D$, we establish
(Theorem \ref{mainresultbanach}) the stretched-exponential bound
\begin{equation}\label{dineq}
|\lambda_n(\l)| <
\frac{W}{r^d} n^{1/2}\ 
r^{\frac{d}{d+1}(d!)^{1/d}n^{1/d}}\quad\text{for all }n\ge1\,,
\end{equation}
where $W:= \sup_{z\in D} \sum_{i} |w_i(z)|$.

We go on to investigate properties of transfer operators acting on
other spaces of holomorphic functions, and prove
(Theorem \ref{ruelleevalueseqthm})
 that the Ruelle
eigenvalue sequence is in a sense \emph{universal}:
for a wide range of domains $D$, and a broad class
of spaces $A(D)$ of holomorphic functions on $D$, the eigenvalue
sequence of $\l:A(D)\to A(D)$ is precisely the Ruelle
eigenvalue sequence.
This universality suggests the 
possibility
of sharpening the estimate (\ref{dineq}), by adapting the proof of
Theorem \ref{mainresultbanach} to some other space $A(D)$.
In particular, the choice of $A(D)$ as the Hardy space
$H^2(D)$ is known to yield a concrete eigenvalue bound for
 $\l:H^2(D)\to H^2(D)$ 
(see \cite{bandtlowjenkinsoncmp}).
Intriguingly, this bound turns out to be \emph{complementary} to
(\ref{dineq}):
in every dimension $d$,
and for every $r<1$,
 (\ref{dineq}) is superior for sufficiently small $n$,
 while the Hardy space bound is superior for sufficiently large $n$.
If $N(r,d)$ denotes the integer such that
 (\ref{dineq}) gives the sharper bound on $|\lambda_n(\l)|$ 
 precisely for $1\le n\le N(r,d)$,
 then both $r\mapsto N(r,d)$ and $d\mapsto N(r,d)$ are increasing
 (cf.~Corollary \ref{combined}, Remark \ref{explicitbounds});
in other words, (\ref{dineq}) is more useful
if the $T_i$ are weakly contracting,
or if the ambient dimension is high.

\section{Transfer operators on favourable spaces of holomorphic functions}
\label{favourablesection}

\begin{notation}
Let $\N$ denote the set of strictly positive integers, and set
$\N_0:=\N\cup\{0\}$.
For $d\in\N$, equip $\C^d$ with the Euclidean 
inner product $(\cdot,\cdot)_{\C^d}$, 
the corresponding
norm $\|\cdot\|_{\C^d}$, and the
induced Euclidean metric, denoted $\delta$. 
For $X\subset \C^d$ we use 
$\Delta_\varepsilon(X)=\all{z\in\C^d}{\delta(z,X)<\varepsilon}$ for the
Euclidean $\varepsilon$-neighbourhood of $X$. 
The set of all bounded domains
(non-empty connected open subsets)
in $\C^d$ will be denoted by $\D_d$. 
For two bounded open sets $\Delta_1,
\Delta_2\subset \C^d$ we write $\Delta_1 \cc \Delta_2$ to mean that
$\overline{\Delta}_1\subset \Delta_2$.
  
Let $B=(B,\norm{\cdot}{B})$ be a Banach space. We often write
$\norm{\cdot}{}$ instead of $\norm{\cdot}{B}$ whenever this does not
lead to confusion. For $X\subset \C^d$ compact and $D\in\D_d$
define 
\begin{align*}
Hol(D,B)&:=\all{f:D\to B}{f \text{ holomorphic}}\\
C(X,B)&:=\all{f:X\to B}{f \text{ continuous}},\,
\norm{f}{C(X,B)}:=\sup_{x\in X}\norm{f(x)}{B} \\
U(D,B)&:=\all{f:\overline{D}\to B}{f\in C(\overline{D},B)\cap
  Hol(D,B)},\,
 \norm{f}{U(D,B)}:=\sup_{z\in \overline{D}}\norm{f(z)}{B}\,.
\end{align*}
Note that $C(X,B)$ and $U(D,B)$ are Banach spaces when equipped with the
indicated norms, while $Hol(D,B)$ is a Fr\'{e}chet space when equipped
with the topology of uniform convergence on compact subsets of $D$.   
If $(B,\|\cdot\|)=(\C,|\cdot|)$ then we
use $C(X)$, $Hol(D)$, and $U(D)$ to denote
$C(X,\C)$, $Hol(D,\C)$, and $U(D,\C)$ respectively.

We use $L(B)$ 
to denote the space of bounded linear operators from a Banach
space $(B,\|\cdot\|)$ to itself, always equipped with the induced operator
norm. 

If $T$ is holomorphic on some $D\in\D_d$,
its derivative at $z\in D$ is denoted by $T'(z)$.
\end{notation}

\begin{defn}
Let $\I$ be a non-empty
countable set.
For  
$D\in\D_d$, a collection 
  $(T_i)_{i\in\I} = (T_i,D)_{i\in\I}$ of
  holomorphic maps $T_i\in U(D,\C^d)$ is called  
  a \textit{holomorphic map system} (on $D$)
  if  
  $\cup_{i\in \I} T_i(D)\cc D$. 
Write
$T_{\underline{i}} :=T_{i_n}\circ\cdots\circ T_{i_1}$
for $\underline{i}=(i_1,\ldots,i_n)\in
\I^n$, $n\in\N$.

For $X\subset \C^d$ compact,
a collection $(T_i)_{i\in \I} = (T_i,X)_{i\in\I}$
of maps $T_i:X\to X$
is a $C^\omega$ \emph{map system} (on $X$) if
there exists $D\in\D_d$ with $X\cc D$
such that each $T_i$ extends holomorphically to $D$ 
and $(T_i,D)_{i\in\I}$ is a holomorphic map system. 
Any such $D$ is called \emph{admissible} for 
the $C^\omega$ map system 
$(T_i,X)_{i\in\I}$.

For $n\in\N$, a $C^\omega$ map system $(T_i,X)_{i\in\I}$
is called \emph{complex $n$-contracting} 
(or simply \textit{complex contracting})
if there exists $D\in\D_d$ with $X\cc D$, such that $T_{\underline
  i}'\in U(D,L(\C^d))$ for every $\underline{i}\in \I^n$ and  
\begin{equation}\label{complexcontracting}
\sup_{\underline i\in\I^n}\norm{T_{\underline
    i}'}{U(D,L(\C^d))}<1\,.
\end{equation}
\end{defn}

Note that if $\I$ is finite then (\ref{complexcontracting}) is implied
by the condition 
$\sup_{\underline i\in\I^n}\|T_{\underline
    i}'\|_{C(X,L(\C^d))}<1$. 

\begin{example}
If $X=[0,1]\subset\C$, define the \emph{Gauss map system}
$(T_i)_{i\in\N}$ by $T_i(x)=1/(i+x)$
(the $T_i$ are the inverse branches to the Gauss map
$x\mapsto 1/x \pmod 1$ on $X$).
This is a $C^\omega$ map system on $X$: for example if 
$D\subset\C$ is the open disc of radius $3/2$ centred at the point $1$
then $(T_i,D)_{i\in\I}$ is a holomorphic map system. 
The system is also complex contracting, because
$\sup_{\underline i\in\I^2}\|T_{\underline i}'\|_{U(D)}=
|T_{(1,1)}'(-1/2)|=4/9<1$
(note we cannot choose $n=1$ in (\ref{complexcontracting}), 
because $T_1'(0)=-1$).
\end{example}

Complex contraction guarantees the existence of an admissible domain,
and this domain may be chosen arbitrarily close to $X$:

\begin{lemma}\label{ccimpliesadmissible}
If a $C^\omega$ map system on $X$ is complex contracting then
there exists a family $\{D_\theta\}_{\theta\in(0,\Theta)}$ 
of admissible domains,
such that $\cap_{\theta\in(0,\Theta)} D_\theta = X$.
\end{lemma}
\begin{proof}
Let $(T_i)_{i\in\I}$ denote the $C^\omega$ map system on $X$.
Choose $n\in\N$ and $D\in\D_d$ such that
$\gamma:=\sup_{\underline i\in\I^n} 
\|T_{\underline i}'\|_{U(D,L(\C^d))}<1$.
From the several variables mean value theorem 
\cite[Thm.~2.3]{avez},
for each
$\underline i\in\I^n$, the map
$T_{\underline i}$ is $\gamma$-Lipschitz,
with respect to Euclidean distance $\delta$,
on any convex
subset of $D$. 
Now set $\beta:=\gamma^{1/n}<1$,
and define the distance 
$$\text{dist}(x,y)= \sup_{\underline i\in\I^{n-1}}
\sum_{k=0}^{n-1} \beta^{n-1-k}\delta(T_{P_k(\underline i)}(x), 
T_{P_k(\underline i)}(y))\,,$$
where for $1\le k\le n-1$, 
$P_k:\I^{n-1}\to \I^k$ denotes the projection 
$P_k{\underline i}=(i_{1}, \ldots, i_{k})$
onto the first $k$ coordinates,
with the convention that 
$ T_{P_0{\underline i}}={\rm id}$.
Note that for each $i\in \I$, the map $T_i$ is
$\beta$-Lipschitz,
with respect to $\text{dist}$,
on any convex
subset of $D$. Moreover, $\rm dist$ and $\delta$ generate the same
topology on $\delta$-compact subsets of $D$. To see this observe that
on the one hand we clearly have 
$\delta(x,y)\leq \beta^{1-n}{\rm dist}(x,y)$ for every $x,y\in
D$. On the other hand, if $K$ is a $\delta$-compact convex subset of
$D$, then by Cauchy's theorem there is $C>0$ such that $\|T_{\underline
  i}'\|_{C(K,L(\C^d))} \le C$ for every 
$\underline{i}\in\I^k$, $1\leq k\leq n-1$. Thus
by the mean value theorem ${\rm dist}(x,y)\leq
C\sum_{k=0}^{n-1}\beta^{n-1-k}\delta(x,y)$ for every $x,y\in K$. 

Since $X$ is compactly contained in the domain $D$,
there exists $\varepsilon>0$ such that the Euclidean
neighbourhood
$\Delta_\varepsilon(X)$
is contained in $D$.
Setting $\Theta:=\varepsilon\beta^{n-1}$, we see that
$D_\theta:=\all{z\in\C^d}{\text{dist}(z,X)<\theta}\subset\Delta_\varepsilon(X)$
for all $\theta\in(0,\Theta]$,
and that $\cap_{\theta\in(0,\Theta]} D_\theta = X$.
If $z\in D_\theta$, and $x\in
X$ satisfies $\text{dist}(z,X)=\text{dist}(z,x)$, then 
$x,z\in\Delta_\varepsilon(x)$, 
a convex subset of $D$, so
$\text{dist}(T_i(z),T_i(x))\le\beta\,\text{dist}(z,x)$.  Therefore
 $ \text{dist}(T_i(z),X) \le\text{dist}(T_i(z),T_i(x)) 
\le\beta\,\text{dist}(z,x) =\beta\,\text{dist}(z,X)$,
and hence
$\cup_{i\in\I} T_i(D_\theta)\subset D_{\beta\theta}\cc D_\theta$,
so $D_\theta$ is admissible for $(T_i,X)_{i\in\I}$.
\end{proof}

\begin{remark}
\label{admissibleremark}
The above proof shows that if a $C^\omega$ map system on $X$ is
complex
\mbox{1-contracting}
then all sufficiently small
Euclidean $\varepsilon$-neighbourhoods
$\Delta_\varepsilon(X)$ 
are admissible.
This is not the case for 
the Gauss map system $T_i(z)=1/(z+i)$ on $X=[0,1]\subset\C$:
no Euclidean $\varepsilon$-neighbourhood 
is admissible, since $\delta(T_1(-\varepsilon),X)
>\varepsilon$.
\end{remark}

\begin{defn} Let $\I$ be a non-empty countable set. 
A \emph{holomorphic weight system} on $D\in\D_d$ is a
  collection $(w_i)_{i\in\I}=(w_i,D)_{i\in\I}$
of holomorphic functions 
(called \emph{weight functions})
$w_i\in U(D)$
such that $\sum_{i\in \I}\norm{w_i}{U(D)}
  <\infty$.

For $X\subset \C^d$ compact,
a collection $(w_i)_{i\in \I}=(w_i,X)_{i\in\I}$
of maps $w_i:X\to \C$
is a \emph{$C^\omega$ weight system} (on $X$) if
there exists $D\in\D_d$ with $X\cc D$
such that $(w_i,D)_{i\in\I}$ 
is a holomorphic weight system.
Any such $D$ is called \emph{admissible}
for $(w_i,X)_{i\in\I}$.

If $(T_i)_{i\in\I}$ is a holomorphic (respectively, $C^\omega$)
map system and
  $(w_i)_{i\in\I}$ is a holomorphic (respectively, $C^\omega$)
weight system then
  $(T_i,w_i)_{i\in\I}$ is called
  a \emph{holomorphic (respectively, $C^\omega$) map-weight system}.
A domain $D\in\D_d$ is called
\emph{admissible} for a $C^\omega$ map-weight system
$(T_i,w_i)_{i\in\I}$ if it is admissible for both $(T_i,X)_{i\in\I}$
and $(w_i,X)_{i\in\I}$.

\end{defn}

With each holomorphic map-weight
system $(T_i,w_i)_{i\in\I}$ we associate a linear operator, 
\begin{equation}
\label{tropdefn}
\l f =\sum_{i\in\mathcal I}w_i\cdot f\circ T_i\,, 
\end{equation}
called the \emph{transfer operator}. 
It will be seen that the transfer operator $\l$ preserves, and acts
compactly upon, 
the following class of spaces 
 of holomorphic functions.

\begin{defn}
\label{favourabledefn}
For $D\in\D_d$, a Banach space $A=A(D)$ of
functions $f:D\to \C$ holomorphic on $D$ is called a
\emph{favourable space of holomorphic functions (on $D$)} if 
\item[\, (i)] for each $z\in D$, the point evaluation functional
$f\mapsto f(z)$ is continuous on $A$, and
\item[\, (ii)]
$A$
contains $U(D)$, and 
the natural embedding\footnote{The embedding $U(D)\hookrightarrow A$ 
is automatically continuous: continuity of point evaluation on both $A$
and $U(D)$ implies that it has closed graph; cf.~the 
proof of Lemma~\ref{continuousembeddinglemma}.}
  $U(D)\hookrightarrow A$ has norm
$1$.
\end{defn}

\begin{remark}
\label{favourablespaceex}
Let $D\in\D_d$. Then $U(D)$ is trivially a favourable space of
holomorphic functions. Other examples include, for $p\in[1,\infty]$, 
\textit{Bergman
  spaces} $L^p_{Hol}(D)$ 
(see \cite[Ch.~I, Cor.~1.7, 1.10]{rangebook}) 
and \textit{Hardy spaces}
$H^p(D)$ (see \cite[Ch.~8.3]{krantz}). If $p=2$ and
$D$ has $C^2$ boundary, then $H^2(D)$ can be identified with the
$L^2(\partial D,\sigma)$-closure of $U(D)$, where $\sigma$ denotes
$(2d-1)$-dimensional Lebesgue on the boundary $\partial D$, normalised
so that $\sigma(\partial D)=1$.
In particular, $H^2(D)$ is a Hilbert space
with inner product given by
$(f,g)=\int_{\partial D}f^*\,\overline{g^*}\,d\sigma$, where, for
$h\in H^2(D)$, the symbol $h^*$ 
denotes the corresponding nontangential limit function 
in $ L^2(\partial D,\sigma)$
--- see \cite[Ch.~1.5 and 8]{krantz}. 
\end{remark}

Recall (see e.g. 
\cite[1.7.1]{pietsch}) 
that a linear operator $L:B\to B$ on a Banach space\footnote{See 
\cite[II, D\'{e}f.~1, p.~3]{grothendieck} 
for the generalisation to locally convex spaces.}
 $B$ is \emph{$p$-nuclear} if there exist sequences
$b_i\in B$ and $l_i\in B^*$ (the strong dual of $B$) with
$\sum_i (\|b_i\|\, \|l_i\|)^p <\infty$,
such that $L(b)=\sum_{i=1}^\infty l_i(b) b_i$ for all $b\in B$.
The operator is \emph{strongly nuclear} 
(or \emph{nuclear of order zero}) if it is $p$-nuclear
for every $p>0$. 
It turns out that certain natural embeddings
between favourable spaces are
strongly nuclear:

\begin{lemma}\label{continuousembeddinglemma}
  Let $D$ and $D'$ be domains in $\C^d$ such that $D'\cc D$. 
Let $A$
  and $A'$ be favourable Banach spaces of holomorphic functions on $D$ and $D'$
  respectively.  Then $A\subset A'$, and the natural embedding
$J:A\hookrightarrow A'$,
defined by $Jf=f|_{D'}$, is strongly nuclear.
\end{lemma} 
\begin{proof}
Choose $D''\in \D_d$ with $D'\cc D''\cc D$, and consider the
natural embeddings
\[ A\overset{J_1}{\hookrightarrow} Hol(D'')
  \overset{J_2}{\hookrightarrow} U(D')
  \overset{J_3}{\hookrightarrow} A'\,.
 \] 
Clearly $J=J_3J_2J_1$. 
The unit ball of $U(D')$ is a neighbourhood in $Hol(D'')$,
so the map $J_2$ is bounded.
But the Fr\'{e}chet space $Hol(D'')$ is nuclear
\cite[II, Cor., p.~56]{grothendieck}, so $J_2$ is
$p$-nuclear for every $p>0$ by 
\cite[II, Cor.~4, p.~39, Cor.~2, p.~61]{grothendieck}. 
It thus suffices to show that $J_1$ 
and $J_3$ are continuous by \cite[I, p.~84, II, p.~9]{grothendieck}.  

Now, $J_3$ is continuous since $A'$ is favourable. 
Finally, to see that $J_1$ is continuous we
note that, by the closed graph theorem 
(see e.g.~\cite[Ch.~III, 2.3]{schaefertvs}), 
it is enough to show that if $f_n\to f$
  in $A$, and $J_1f_n\to g$ in $Hol(D'')$, then $g=J_1f=f|_{D''}$.  Since
  point evaluation is continuous on $A$, $f_n(z)\to f(z)$ for all
  $z\in D$ and in particular for all $z\in D''$.  But point
  evaluation is also continuous on $Hol(D'')$, so $f_n(z)=J_1f_n(z)\to g(z)$
as $n\to\infty$ for all $z\in D''$.  Therefore $g=f|_{D''}$.
\end{proof}

Favourable spaces $A$ are always
invariant under the transfer operator $\l$, and 
the restricted operator (henceforth denoted by $\ls{A}$)
is always compact, indeed strongly nuclear:

\begin{prop}\label{compactness}
  Let $(T_i,w_i,D)_{i\in\I}$ be a holomorphic map-weight system.
The corresponding transfer operator leaves invariant
every favourable space $A$ of holomorphic functions on $D$, and
  $\ls{A}:A\to A$ is strongly nuclear.  
\end{prop}
\begin{proof}
Choose $D'\in\D_d$ with
$
\cup_{i\in\I} T_i(D) \cc D' \cc D
$.
First we observe that $\hat{\l}f:= \sum_{i\in\I} w_i\cdot f\circ T_i$
defines a continuous operator $\hat{\l}:U(D')\to
 U(D)$.
To see this, fix $f\in U(D')$ and note that $w_i\cdot f\circ
 T_i\in U(D)$ with $\norm{w_i\cdot f\circ
 T_i}{U(D)}\leq \norm{w_i}{U(D)}\norm{f}{U(D')}$ 
for every $i\in\I$. But since     
$\|\hat{\l}f\|_{U(D)}\leq
\sum_{i\in\I}\norm{w_i}{U(D)}\norm{f}{U(D')}$ and
$\sum_{i\in\I}\norm{w_i}{U(D)}<\infty$ by hypothesis, we conclude that
$\hat{\l} f\in U(D)$ and 
that $\hat{\l}$ is continuous. 

Since $A$ is favourable,
the embedding $\hat{J}:U(D)\hookrightarrow
A$ is continuous,
and the embedding
$J:A\hookrightarrow U(D')$ is $p$-nuclear for every $p>0$ 
by Lemma \ref{continuousembeddinglemma}. 
Moreover, if $f\in A$ then $\l f=\hat{J}\hat{\l}Jf\in A$. Thus 
$A$ is $\l$-invariant,
and the operator 
$\ls{A}=\hat{J}\hat{\l}J$
is $p$-nuclear for any $p>0$. 
\end{proof}

\begin{remark}
\label{compactnessremark}
Strong nuclearity of the transfer operator on spaces of holomorphic functions
is not new (the original
result of this kind is
\cite{ruelleinventiones}, but see also 
e.g.~\cite{glz,jpadvmath, mayergreenbook}); the novelty of Proposition
\ref{compactness} is in the breadth of spaces covered.\footnote{Actually
the result
can be further extended to certain locally convex spaces of 
holomorphic functions, including $Hol(D)$.}
\end{remark}

\section{Eigenvalue bounds}
\label{weylsection}

For favourable $A$, the compactness of $\ls{A}$ means its spectrum
consists of
a countable set of
eigenvalues, each with
finite algebraic multiplicity, together 
with a possible accumulation point at 0.
We wish to obtain bounds on the \emph{eigenvalue sequence}
$\lambda(\ls{A}):=\set{\lambda_n(\ls{A})}_{n=1}^\infty$,
i.e.~the sequence of all
eigenvalues of $\ls{A}$ counting algebraic
multiplicities and ordered by decreasing modulus.\footnote{By 
convention
distinct eigenvalues with the same modulus can be written in 
any order (see e.g.~\cite[3.2.20]{pietsch}).}

If $L:B_1\to B_2$ is a continuous operator
between Banach spaces then for $k\ge1$, its \emph{$k$-th approximation
  number} $a_k(L)$ is defined as
$$
a_k(L)= \inf \all{\norm{L-K}{}}{K:B_1\to B_2\text{ linear with }{\rm
    rank}(K)<k}\,.
$$

\begin{prop}
\label{weylforLbanach}
For a $C^\omega$ map-weight system 
$(T_i,w_i)_{i\in\I}$ such that $(T_i)_{i\in\I}$ is complex 
contracting, and a favourable space $A=A(D)$ such that $D\in \D_d$
is admissible,
\begin{equation}
\label{weylbanachforL}
|\lambda_n(\ls{A})| \le W n^{1/2}\prod_{k=1}^n a_k(J)^{1/n}
\quad\text{for all }n\ge1\,,
\end{equation}
where 
$W:=\sup_{z\in D}\sum_{i\in\I} |w_i(z)|$,
$D'\in\D_d$ is such that
$
 \cup_{i\in\I}T_i(D) \subset D' \cc D 
$, and
$J:A(D)\hookrightarrow U(D')$ is
the canonical embedding.
\end{prop}
\begin{proof}
Since $A(D)$ is favourable, the embedding
$\hat{J}:U(D)\hookrightarrow A(D)$ is continuous of norm $1$.
Observe that 
$\hat{\l}f= \sum_{i\in\I}w_i\cdot
 f\circ T_i$ defines a continuous operator $\hat{\l}:U(D')\to
 U(D)$ (see the proof of Proposition~\ref{compactness}) with  
$\|\hat{\l}\| \le W$. To see the latter note that for $f\in U(D')$ we
 have $|f(T_i(z))|\leq \norm{f}{U(D')}$ for every $z\in D$,
 $i\in\I$; thus by the maximum principle 
$\|\hat{\l}f\|_{U(D)}=\sup_{z\in D}|(\hat{\l}f)(z)|\leq \sup_{z\in D}
  \sum_{i\in\I}\abs{w_i(z)}\,\abs{f(T_i(z))}\leq
 W\norm{f}{U(D')}$. 

Now clearly 
$\ls{A}=\hat{J}\hat{\l}J$, so
\begin{equation}
\label{embeddinglift}
a_k(\ls{A})\le \|\hat{J}\hat{\l}\|a_k(J) \le Wa_k(J)
\quad\text{for all }k\ge1\,,
\end{equation}
since in general $a_k(L_1L_2)\le \norm{L_1}{} a_k(L_2)$ whenever $L_1$
and $L_2$ are bounded operators between Banach spaces (see
\cite[2.2]{pietsch}).
Moreover, since $\ls{A}$ is
compact,
Weyl's inequality
 (see e.g.~\cite{hinrichs})
asserts that
$
 \prod_{k=1}^n \left| \lambda_k(\ls{A})\right|\, \leq \, n^{n/2}\,
\prod_{k=1}^n a_k(\ls{A})
$ for every $n\in\N$.\footnote{This is a 
Banach space version of Weyl's original inequality
in Hilbert space (see \cite{weyl}). Note that the constant $n^{n/2}$ is optimal
(see \cite{hinrichs}).} 
Together with (\ref{embeddinglift}) this yields
(\ref{weylbanachforL}),
because 
$|\lambda_n(\ls{A})|\le \prod_{k=1}^n |\lambda_k(\ls{A})|^{1/n}$.
\end{proof}

Taking $A(D)=U(D)$, the \emph{Ruelle eigenvalue sequence}
$\lambda(\ls{U(D)})$ can be bounded as follows:

\begin{theorem}
\label{mainresultbanach}
  Suppose the Euclidean ball $D\subset\C^d$
is an admissible domain for a $C^\omega$ map-weight system 
$(T_i,w_i)_{i\in\I}$,
and that $\cup_{i\in\I}T_i(D)$
is contained in the concentric ball whose radius is $r<1$ times 
that of $D$.
Setting $W:=\sup_{z\in B}\sum_{i\in\I}|w_i(z)|$, the 
Ruelle eigenvalue sequence 
$\lambda(\ls{U(D)})$ can be bounded by
\begin{equation}
\label{generalevaluebanach}
|\lambda_n(\ls{U(D)})| <
\frac{W}{r^d} n^{1/2}\ 
r^{\frac{d}{d+1}(d!)^{1/d}n^{1/d}}\quad\text{for all }n\ge1\,.
\end{equation}
If $d=1$ then
\begin{equation}
\label{d=1evaluebanach}
 |\lambda_n(\ls{U(D)})|\leq 
W n^{1/2}\  r^{(n-1)/2}\quad\text{for all }n\ge1\,.
\end{equation}
\end{theorem}
\begin{proof}
Without loss of generality let $D=D_1$
be the open unit ball,
and let the smaller concentric
ball be $D_r$, the ball of radius $r$ centred at $0$.
Let $J:U(D_1)\hookrightarrow U(D_r)$ 
be the canonical embedding. 
From \cite[Prop.~2.1 (a)]{farkov2} it follows that  
$a_l(J)\leq r^{\bt_l}$,
where $\bt_l:=k$ for
$\binom{k-1+d}{d}<l\le \binom{k+d}{d}$,
hence 
$\prod_{l=1}^n a_l(J)^{1/n} \le r^{\frac{1}{n}\sum_{l=1}^n\bt_l}$.
If $d=1$ then
$\frac{1}{n} \sum_{l=1}^n \bt_l
=\frac{1}{n} \sum_{l=1}^n (l-1) = (n-1)/2$,
and (\ref{d=1evaluebanach})
 follows from
(\ref{weylbanachforL}).
More generally
$
\bt_l \ge (d!)^{1/d} l^{1/d} -d
$,
so that
\begin{equation*}
\label{betanbound}
\frac{1}{n} \sum_{l=1}^n \bt_l
\ge -d +  (d!)^{1/d} \frac{1}{n} \sum_{l=1}^n l^{1/d} 
>  -d +  (d!)^{1/d} \frac{d}{d+1}n^{1/d}
\end{equation*}
using the estimate
$\sum_{l=1}^n l^{1/d} >  \int_{x=0}^n x^{1/d} =\frac{d}{d+1}
n^{1+1/d}$,
and (\ref{generalevaluebanach})
 follows from
(\ref{weylbanachforL}).
\end{proof}

\section{Universality of the Ruelle eigenvalue sequence}

If $(T_i,w_i)_{i\in\I}$ is a $C^\omega$ map-weight system with complex
contracting $(T_i)_{i\in \I}$ then, in view of 
Lemma~\ref{ccimpliesadmissible} and
Proposition~\ref{compactness}, there is some freedom in the choice of
an admissible $D$, and a favourable space $A=A(D)$ on which to consider
the transfer operator $\ls{A}$. 
The purpose of this section is to show
that the eigenvalue sequence of $\ls{A}$ is 
in fact independent of $A$:
it is always equal to the Ruelle eigenvalue sequence 
$\lambda(\ls{U(D)})$ (see Corollary \ref{spectrumindept}).
For this, we first
require some facts from the Fredholm theory
originally developed by Grothendieck \cite{grothendieck}. 
If $B$ is a Banach space, we denote by $N_p(B)$ ($p>0)$ the
quasi-Banach operator ideal of $p$-nuclear operators  on $B$
(cf.~\cite[D.1.4, 1.7.1]{pietsch}). 
If $p\leq 2/3$ then $N_p(B)$ admits a unique continuous trace $\tau$ and
a unique continuous determinant $\det$ (see \cite[1.7.13, 4.7.8,
4.7.11]{pietsch}), related for a fixed $L\in
N_p(B)$ by  
\begin{equation}
\label{tracedetrelation}
\det(I-zL)=\exp\left (-\sum_{n=1}^\infty \frac{z^n}{n}\tau(L^n)\right
)\,,
\end{equation} 
for all $z\in\C$ in a suitable neighbourhood of $0$ 
(see \cite[4.6.2]{pietsch}). Moreover, both
$\tau$ and $\det$ are spectral, which means that
$\tau(L)=\sum_{n=1}^\infty \lambda_n(L)$ and that, counting multiplicities,
the zeros of the entire
function $z\mapsto \det(I-zL)$ are precisely the reciprocals of
the eigenvalues of $L$ (see \cite[4.7.14, 4.7.15]{pietsch}).

\begin{defn}
To any holomorphic map-weight system
$(T_i,w_i)_{i\in\I}$, the associated 
\emph{dynamical determinant}
is the entire function $\Delta:\C\to\C$, defined  
for all $z$ of sufficiently small modulus by
\begin{equation}\label{ruelledet}
\Delta(z)=
\exp
\left (-\sum_{n\in\N}\frac{z^n}{n}\sum_{\underline{i}\in \I^n} 
\frac{w_{\underline{i}}
(z_{\underline{i}})}{\det(I-T_{\underline{i}}'(z_{\underline{i}}))}\right
)\,,
\end{equation}
where
$w_{\underline i}:=\prod_{k=1}^nw_{i_k}\circ T_{P_{k-1}{\underline
    i}}$,
$P_k:\I^n\to \I^k$ denotes the projection 
$P_k{\underline i}=(i_1, \ldots, i_k)$
with the convention that 
$ T_{P_0{\underline i}}={\rm id}$,
and $z_{\underline i}$ denotes the 
(unique, by
\cite{earlehamilton}) fixed-point
of $T_{\underline i}$ in $D$.
\end{defn}

Ruelle \cite{ruelleinventiones} showed
that $\Delta$ is the determinant
of the strongly nuclear operator $\l:U(D)\to U(D)$.
Therefore,
if the zeros $z_1,z_2,\ldots$
of $\Delta$ are listed according to increasing modulus and
counting multiplicity, then the reciprocal sequence
$\{z_n^{-1}\}_{n=1}^\infty$
is precisely the Ruelle eigenvalue sequence.

\begin{theorem}\label{ruelleevalueseqthm}
Let $(T_i,w_i,D)_{i\in\I}$ be a holomorphic map-weight system.
Then the associated transfer operator preserves every favourable space
of holomorphic functions on $D$,
and its 
determinant on each of these spaces
is precisely the dynamical determinant $\Delta$.
\end{theorem}
\begin{proof}
Comparison of 
(\ref{tracedetrelation})
and (\ref{ruelledet})
means 
we require the trace formula\footnote{This trace formula (\ref{traceformula}) generalises
the original one of Ruelle \cite{ruelleinventiones}
for $A=U(D)$, as well as that of 
Mayer \cite{mayersmf, mayercompop, mayergreenbook}.
Our method of proof
is rather direct,
reducing to a simple Hilbert space computation;
in particular, we do not need to
explicitly evaluate the eigenvalues
of each weighted composition operator $f\mapsto w_{\underline{i}}\cdot
f\circ T_{\underline{i}}$ 
(a more complicated procedure, particularly
in higher dimensions,
cf.~\cite[\S III]{mayercompop}).}
\begin{equation}\label{traceformula}
\tau(\l_A^n)=
\sum_{\underline{i}\in \I^n} 
\frac{w_{\underline{i}}
(z_{\underline{i}})}{\det(I-T_{\underline{i}}'(z_{\underline{i}}))}
\quad\text{for all }n\ge1\,,
\end{equation}
for every favourable space $A$ on the admissible
domain $D$.

First consider the
holomorphic map-weight system $(T,w,D)$ consisting of a single map and
weight.
Since $T(D)\cc D$, the Earle-Hamilton theorem \cite{earlehamilton}
implies that $T$ has a unique
fixed-point $z_0\in D$, and the
eigenvalues of $T'(z_0)$ lie in the open unit disc 
\cite[Thm.~1]{mayercompop}. 
If $\l f=w\cdot f\circ T$ is the
corresponding transfer operator, 
we claim that
\begin{equation}\label{simplesingletrace}
\tau(\ls{A})=\frac{w(z_0)}{\det(I-T'(z_0))}\,.
\end{equation}
The admissibility of $D$ and favourability
of $A=A(D)$ are invariant under affine coordinate
changes, and $\tau$ is
invariant under continuous
similarities, so we may assume that $z_0=0$
and $\norm{T'(0)}{L(\C^d)}<1$.  Therefore, by
Lemma~\ref{ccimpliesadmissible} and Remark~\ref{admissibleremark},
there exists $R>0$ such that, for $r\in (0,R)$,
the radius-$r$ Euclidean ball $B_r$ centred
at $0$ is admissible. 

Let $H^2_r=H^2(B_r)$ denote the Hardy space on
$B_r$,
a favourable Hilbert space 
(see Remark \ref{favourablespaceex})
with inner product 
$(f,g)_{H^2_r}=\int_{S_r}f^*\, \overline{g^*}\,d\sigma_r$, 
where $S_r=\partial B_r$,
$\sigma_r(S_r)=1$,
and with orthonormal basis
(cf.~\cite[Prop.~1.4.8, 1.4.9]{rudin})
$\all{p_{\underline{n},r}}{\underline{n}\in \N_0^d}$,
where $p_{\underline{n},r}(z)=K_{\underline{n}}r^{-|\underline{n}|}z^{\underline{n}}$
and 
$ 
K_{\underline{n}}
= \sqrt{\frac{(|\underline{n}|+d-1)!}{(d-1)!\,\underline{n}!}}$,  
$\un=(n_1,\ldots,n_d)$, $z^{\un}=z_1^{n_1}\cdots z_d^{n_d}$,
 $\un!=n_1!\cdots n_d!$, $|\un|=n_1+\cdots +n_d$.

The canonical 
embedding $J_r:A\hookrightarrow H^2_r$
has dense range, because
complex polynomials are dense in $H^2_r$, 
and $J_r\ls{A}=\ls{H^2_r}J_r$.
An intertwining argument of Grabiner \cite[Lem.~2.3]{grabiner} 
then implies that
$\lambda(\ls{A})=\lambda(\ls{H^2_r})$,
and hence that
$\tau(\ls{A})=\tau(\ls{H^2_r})$ because $\tau$ is spectral. 
The strong nuclearity of $\ls{H^2_r}$ 
means it is trace class,
so $\tau(\ls{H^2_r})$ equals the sum of the diagonal
    entries of the matrix representation of $\ls{H^2_r}$ with respect
      to an orthonormal basis. Thus, for any $r\in(0,R)$, 
\begin{multline*}
\tau(\ls{A})=\tau(\ls{H^2_r})=\sum_{\underline{n}\in\N_0^d}(\l
p_{\underline{n},r},
p_{\underline{n},r})_{H^2_r}=
\int_{S_r}w(z)\sum_{\underline{n}\in\N_0^d}K_{\underline{n}}^2r^{-2|\underline{n}|} T(z)^{\underline{n}}\,\overline{z}^{\underline{n}}\,d\sigma_r(z) \\
=\int_{S_r} \frac{w(z)}{\left (1-(r^{-1}T(z), r^{-1}z)_{\C^d}\right )^d}\,d\sigma_r(z)
= \int_{S_1} \frac{w(rz)}{\left (1-(r^{-1}T(rz), z)_{\C^d}\right )^d}\,d\sigma_1(z)\,.
\end{multline*}
Letting
$r\to 0$ gives
\[ \tau(\ls{A})=\int_{S_1} \frac{w(0)}{\left (1-(T'(0)z,
    z)_{\C^d}\right )^d}\,d\sigma_1(z)=\frac{w(0)}{\det(I-T'(0))} \]
by an elementary integration,
and (\ref{simplesingletrace}) is proved.

Returning to the case of the holomorphic map-weight system
$(T_i,w_i)_{i\in\I}$,
the factorisation argument used in the proof of
Proposition~\ref{compactness} shows that for $n\in \N$, the series
$\sum_{\underline{i}\in \I^n} \l_{\underline{i}}$
  converges in $N_{2/3}(A)$
to $\ls{A}^n$, where $\l_{\underline{i}}:A\to A$ is given by 
$\l_{\underline{i}}f=w_{\underline{i}}\cdot f\circ T_{\underline{i}}$.
Since $\tau$ is continuous, 
$\tau(\ls{A}^n)=\sum_{\underline{i}\in\I^n}\tau(\l_{\underline{i}})$, 
and the required
trace formula (\ref{traceformula})
follows from 
(\ref{simplesingletrace}).
\end{proof}

\begin{cor}
\label{spectrumindept}
Let $(T_i,w_i)_{i\in\I}$ be a $C^\omega$ map-weight system  
such that $(T_i)_{i\in\I}$ is complex contracting. 
Then the associated transfer operator preserves every favourable space on
every admissible domain, and its eigenvalue sequence 
on each of these spaces is precisely the Ruelle eigenvalue sequence.
\end{cor}

In view of Corollary~\ref{spectrumindept}, 
the Ruelle eigenvalue sequence
associated with a complex contracting
$C^\omega$ map-weight system 
will
henceforth be denoted 
simply by $\lambda(\l)=\set{\lambda_n(\l)}_{n=1}^\infty$. 

\begin{cor}\label{combined}
Under the hypotheses of Theorem \ref{mainresultbanach},
the 
Ruelle eigenvalue sequence 
$\lambda(\l)$ can be bounded by
\begin{equation}\label{combinedeq}
|\lambda_n(\l)| <
\min \left( n^{1/2}\ ,\
\frac{\sqrt{d}}{(1-r^2)^{d/2}}\ n^{(d-1)/(2d)} \right) 
\frac{W}{r^d} 
r^{\frac{d}{d+1}(d!)^{1/d}n^{1/d}}\,.
\end{equation}
\end{cor}
\begin{proof}
Hardy space $H^2(D)$ is favourable,
so Corollary~\ref{spectrumindept} implies that
$\lambda(\ls{H^2(D)})$ is the Ruelle eigenvalue sequence.
The bound
$$
|\lambda_n(\l)| <
\frac{W \sqrt{d}}{r^d (1-r^2)^{d/2}}\ n^{(d-1)/(2d)}\ 
r^{\frac{d}{d+1}(d!)^{1/d}n^{1/d}}
$$ 
then follows from \cite[Thm.~1]{bandtlowjenkinsoncmp}.
The other part of (\ref{combinedeq})
is immediate from
Theorem \ref{mainresultbanach}.
\end{proof}

\begin{remark}\mbox{}
\label{explicitbounds}
For a given $(T_i,w_i)_{i\in\I}$, 
if $r<1$ is chosen as small as possible
then the part of (\ref{combinedeq})
arising from \cite{bandtlowjenkinsoncmp}
is asymptotically superior as $n\to\infty$.
For sufficiently small $n$,
the part of (\ref{combinedeq})
arising from
Theorem \ref{mainresultbanach} is sharper.
For example, in dimension $d=1$ this latter
bound on $|\lambda_n(\l)|$
is superior whenever $n^2 < 1/(1-r^2)$; this is always
the case for $n=1$, and may be true for many $n$
if $r$ is large (i.e.~the map system is weakly contracting).
\end{remark}

\begin{ack}
We are grateful to B.~Carl,
H.~K\"onig, V.~Mascioni and A.~Pietsch for advice on the sharpness of
Weyl's inequality in Banach spaces, and K.~Osipenko, A.~Pinkus and
M.~Stessin for bibliographical assistance concerning the approximation
numbers for the embedding of $U(B_1)$ in $U(B_r)$.
Both authors were partially supported by EPSRC grants
GR/R64650/01 (2001--2003) and  GR/S50991/01 (2004--2007), 
and the second author was partially supported by an EPSRC
  Advanced Research Fellowship.
\end{ack}

\end{document}